\documentclass{amsart}

\usepackage{amsmath,amssymb,graphicx}
\usepackage[colorlinks]{hyperref}

\newcommand{\Z}{\mathbb Z}
\newcommand{\C}{\mathcal C}
\newcommand{\f}{\mathfrak f}
\newcommand{\Li}{\mathcal L}

\newtheorem{lemma}{Lemma}[section]
\newtheorem{proposition}[lemma]{Proposition}

\newtheorem{theorem}[lemma]{Theorem}

\begin{document}
\title[Teaching perspectives of the coin problem]{Teaching perspectives of the Frobenius coin problem of two denominators}
\author[G.~Kapetanakis]{Giorgos Kapetanakis}
\address{Department of Mathematics, University of Thessaly, 3rd km Old National Road Lamia-Athens, 35100, Lamia, Greece}
\email{kapetanakis@uth.gr}
\author[I.~Rizos]{Ioannis Rizos}
\address{Department of Mathematics, University of Thessaly, 3rd km Old National Road Lamia-Athens, 35100, Lamia, Greece}
\email{ioarizos@uth.gr}
\date{\today}

\begin{abstract}
Let $a,b$ be positive, relatively prime, integers. Our goal is to characterize, in an elementary way, all positive integers $c$  that can be expressed as a linear combination of $a,b$ with non-negative integer coefficients and discuss the teaching perspectives of our methods.
\end{abstract}



\maketitle
\section{Introduction}
In order to contribute to the creation of didactical situations that could increase the active involvement of undergraduate students in the course of Number Theory and at the same time offer them the possibility to develop their own solving strategies, we study the well-known \emph{Frobenius coin problem} for two denominators. Our experience suggests that problems like the one above, which can be reduced to solving a Diophantine equation and are essentially open-ended problems (in the sense that may have several paths to resolution), can really initiate the feeling of engagement in the participants and help them to upgrade their problem-solving skills \cite{rizosgkrekas23}.

In parallel the mathematical elaboration aims into guiding the students to discovering known proofs of the problem themselves, using various approaches, with encouragement and strategic guidance from the instructor. 
\section{Motivation and notation}\label{motivation}
Number Theory is a mathematical tradition with a long history and its own peculiar and always attractive problems. The peculiarity of the problems of Number Theory consists in a rather strange character; while these problems appear most of the time as understandable in their formulation by a high school student, their solution often requires knowledge of so-called ``advanced mathematics'' -or at least of algebra taught in the first years of mathematics departments- and sometimes may even be unattainable (except in special cases) by the mathematical community. At the same time, the same problems are usually empirically explorable by simple arithmetic operations, done from memory or with the help of a simple computer program. This makes it easier to explore special cases and to refute any incorrect conjectures or, conversely, to strengthen the belief in the truth of some carefully formulated claims following the negations, following Lakatos’ heuristic epistemological model for mathematics \cite{lakatos76}.
 
The last possibility constitutes, in our opinion, the main epistemological advantage of Number Theory over the other branches of elementary mathematics, in terms of introducing students to the actual proof process. For, as Popper~\cite{popper02} and Lakatos~\cite{lakatos76} have demonstrated, in the Philosophy of Science \emph{the logic of scientific discovery} cannot be separated from \emph{the logic of proof or justification}. More specifically, in order to determine whether students are able to solve simple word problems that require basic knowledge of Number Theory and Algebra, entering into an elementary proof process, we suggest the following tasks:

\subsection*{Task 1} The well-known \emph{Frobenius coin problem} for two denominators.
\begin{enumerate}
    \item Which is the largest monetary amount that cannot be obtained using only coins of 3 and 7 units?
    \item Which is the largest monetary amount that cannot be obtained using only coins of 12 and 25 units?
    \item Taking into consideration the two cases above, can you guess what the largest monetary amount that cannot be obtained using only coins of $a$ and $b$ units, where $\gcd(a, b) = 1$, might be?
    \item Prove by mathematical induction (or by any other method you consider appropriate) that the largest monetary amount that cannot be obtained using only coins of $a$ and $b$ units, where $\gcd(a, b) = 1$, is $ab - (a + b)$.
\end{enumerate}

\subsection*{Task 2} Given the equation 
\begin{equation}\label{1}
5x+8y=d,
\end{equation}
with unknowns $x,y$ and parameter $d\in\Z_{\geq 0}$. Find the maximum value of the parameter $d$ for which Eq.~\eqref{1} is impossible in $\Z_{\geq 0}$, i.e., there is no pair of non-negative integers $x,y$ that verifies it.

\subsection*{Task 3} With the help of dynamic geometry software (e.g. {\sc GeoGebra}) draw in an orthonormal coordinate system the (variable) line
\[ \varepsilon : 5x+8y=d , \]
where the parameter $d$ takes positive integral values less than or equal to $40$. Find the maximum value of the parameter $p$ for which $\varepsilon$ does not pass through a point with coordinates of natural numbers.
%
%
%
%

The study of the coin problem of two denominators translates into studying the set
\[ \C_{a,b} = \{ d\in\Z_{\geq 0} : d = ax+by\text{ for some } x,y\in\Z_{\geq 0} \} \]
and identifying the largest integer that is not in this set. This integer is known as the \emph{Frobenius number} and we denote it by $\f_{a,b}$. Regarding this number, it is well-known that 
\begin{equation}\label{f}
\f_{a,b} = ab-(a+b).
\end{equation}
For a more comprehensive presentation of the coin problem we refer the interested reader to the excellent textbook of Alfons\'in~\cite{alfonsin05} and the references therein. Also, we note that the problem of interest can be generalized by admitting an arbitrary number of denominators in an obvious way. However, the complexity of the resulting problem increases, in such a way, that even the computing $\f_{a,b,c}$ (for three denominators) is extremely complicated and it widely accepted that a closed formula for this number is out of reach, see \cite{alfonsin05}.
So, due to the large complexity of the case of more than two denominators, we will restrict ourselves to the simplest case of two denominators.

To be more precise, we suggest methods, based on known proofs, whose goal is to guide students into discovering the existence of $\f_{a,b}$, into computing $\f_{a,b}$ and into characterizing $\C_{a,b}$, instead of immediately provide  a formal proof. Towards this end, we will use three different approaches, one based on elementary number theory, one based on geometry and one using a newly discovered inductive method.
\section{Mathematical elaboration}
It is clear that we will be working with the linear Diophantine equation
\begin{equation}\label{dioph}
d=ax+by,
\end{equation}
where, we assume that $0<a,b$ and $\gcd(a,b)=1$. We shall also assume that $d\geq 0$, since we are interested in finding positive integer solutions of Eq.~\eqref{dioph} and it is clear that the case $d<0$ is impossible. We will first focus on the number theoretic approach and then on the geometric one and we note that these approaches are inspired by the first and third proof of \cite[Theorem~2.1.1]{alfonsin05} respectively. Last but not least, we will attempt to heuristically communicate a recently discovered inductive proof \cite{kapetanakisrizos23b}. 
We note that \cite{alfonsin05} contains four proofs of this theorem. We chose to only present the first and the third one, as the second seems to contain a mistake\footnote{In particular, it is based on the statement that there exist positive integers $x,y$, such that $ax+by=1$, which is clearly not true.} and the fourth one incorporates analytic concepts whose direct relation with the problem may be clear to university freshmen.

In all cases, we assume that the students are aware of the following celebrated result from Number Theory:
\begin{theorem}[B\'ezout's identity] \label{bezout}
  Let $a,b\in\Z$ with $c=\gcd(a,b)$. There exist some $x',y'\in\Z$, such that $ax' +by'  = c$.
\end{theorem}
We note that one such pair $(x',y')$ can be identified with the Euclidean algorithm. Also, it is obvious that Theorem~\ref{bezout} implies that Eq.~\eqref{dioph} is solvable. In particular, Theorem~\ref{bezout} implies that $ax' +by'  = 1$ for some $x',y'$, hence, $a(dx')+b(dy') = d$, that is, $(x_0,y_0) = (dx',dy')$ is a solution of Eq.~\eqref{dioph}.

However it is neither guaranteed that the numbers $x_0,y_0$ are both positive, nor that the solution is unique. In fact, the students can easily observe, with the help of easy examples, that the former is sometimes true and the latter is not true at all. Next, we will identify the full set of solutions. Towards this goal, let $(x_1,y_1)$ and $(x_2,y_2)$ be two solutions of Eq.~\eqref{dioph}. We observe that:
\[ \left. \begin{array}{rl}
  d & = ax_1+by_1 , \\
  d& = ax_2+by_2
\end{array} \right\} \Rightarrow a(x_1-x_2) = b(y_2-y_1) \stackrel{\gcd(a,b)=1}{\Longrightarrow} a \mid y_2-y_1 \Rightarrow y_2=y_1+ka , \]
for some $k\in\Z$. It follows that $x_2=x_1-kb$, thus, $(x_2,y_2) = (x_1-kb,y_1+ka)$. Further, it is not hard to check, that for every $k\in\Z$, the pair $(x-kb,y+ka)$ is a solution of Eq.~\eqref{dioph}. We have established the following:
\begin{proposition} \label{gen_sol}
 Let $a,b\in\Z$, not both zero, with $\gcd(a,b)=1$. If $x_0,y_0\in\Z$ are such that $ax_0 +by_0  = d$, then the set of solutions of Eq.~\eqref{dioph} is the (infinite) set
\[ S = \{ (x_0-kb,y_0+ka) : k\in\Z\} . \]
\end{proposition}
\subsection{Number theoretic approach}
A couple of key observations that derive from Proposition~\ref{gen_sol} are the following:
\begin{enumerate}
  \item It is clear that Eq.~\eqref{dioph} admits a non-negative solution, if $x$ is non-negative, but as small as possible and $y$ turns out to be also non-negative.
  \item If $(x,y)$ is a solution of Eq.~\eqref{dioph}, then $ax+by = d \Rightarrow x\equiv d\cdot a^{-1} \pmod b$, that is, the smallest non-negative value for $x$ is the number $x_{(d)}$ that satisfies $0\leq x_{(d)} < b$ and $x_{(d)}\equiv d\cdot a^{-1} \pmod b$. Further, observe that, since $d$ takes arbitrary values, $x_{(d)}$ can in fact take any value within the described range.
\end{enumerate}
From these observations, we get that the worst case scenario is for $(x_{(d)},y_{(d)}) = (b-1,-1)$ to be a solution and this occurs exactly when $d=ab-a-b$. Eq.~\eqref{f} follows.

In addition, we obtain an easy characterization of $\C_{a,b}$. In particular, given $a,b\in\Z_{>0}$, with $\gcd(a,b)=1$ and some $0\leq d< ab-a-b$, we can check whether $d\in\C_{a,b}$ as follows:
\begin{enumerate}
  \item Compute $x_{(d)}$ as the unique number such that
   \[ 0\leq x_{(d)}<b\text{ and }x_{(d)}\equiv d\cdot a^{-1}\pmod b. \]
  \item Compute $y_{(d)} = (d-ax_{(d)})/b$.
  \item If $y_{(d)}\geq 0$, then $d\in\C_{a,b}$, otherwise $d\not\in\C_{a,b}$.
\end{enumerate}
\subsection{Geometric approach}
For this approach, we will be working on the Cartesian plane. On the plane, the problem in question, translates into whether the line
\[ \Li_d : d = ax+by \]
meets the lattice
\[ \Lambda := \{ (x,y) : x,y\in\Z \} \]
in the first quadrant or in one of the positive semi-axes. Note that Proposition~\ref{gen_sol} ensures that $\Li_d$ and $\Lambda$ meet at an infinite number of points (for fixed $d$), while it is clear that every point of $\Lambda$ is contained in some line $\Li_d$ (for some $d$). Furthermore, since $a,b,d$ are all assumed to be non-negative, the line $\Li_d$ intersects the first quadrant and the axes at a segment, that is, the number of points in question will be finite. Also, the line $\Li_d$ does not pass through the fourth quadrant, thus not both coordinates of the points of $\Li_d$ can be negative.

\begin{figure}[h]
\begin{center}
  \includegraphics[width=\textwidth]{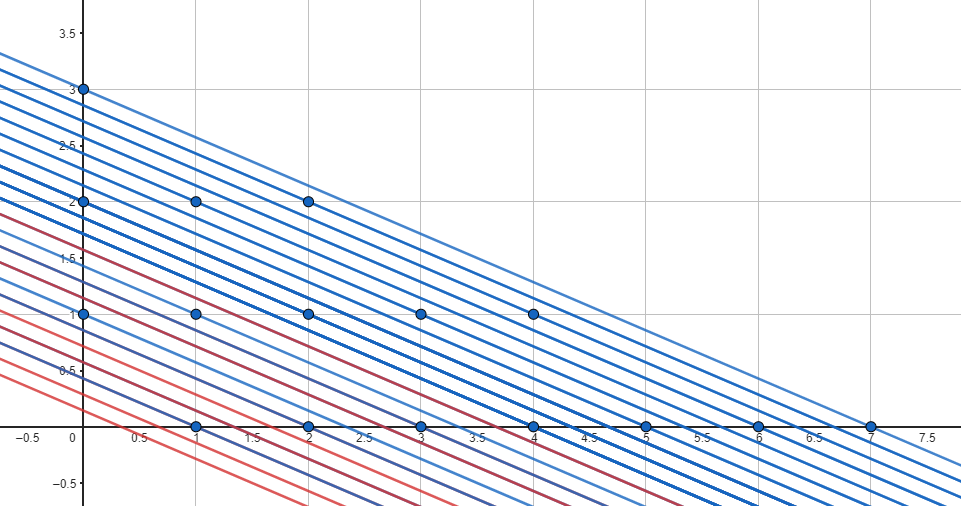}
\end{center}
\caption{{\sc GeoGebra}'s visualization of $3x+7y=d$ for $1\leq d\leq 30$.\label{f1}}
\end{figure}

\begin{figure}[h]
\begin{center}
  \includegraphics[width=\textwidth]{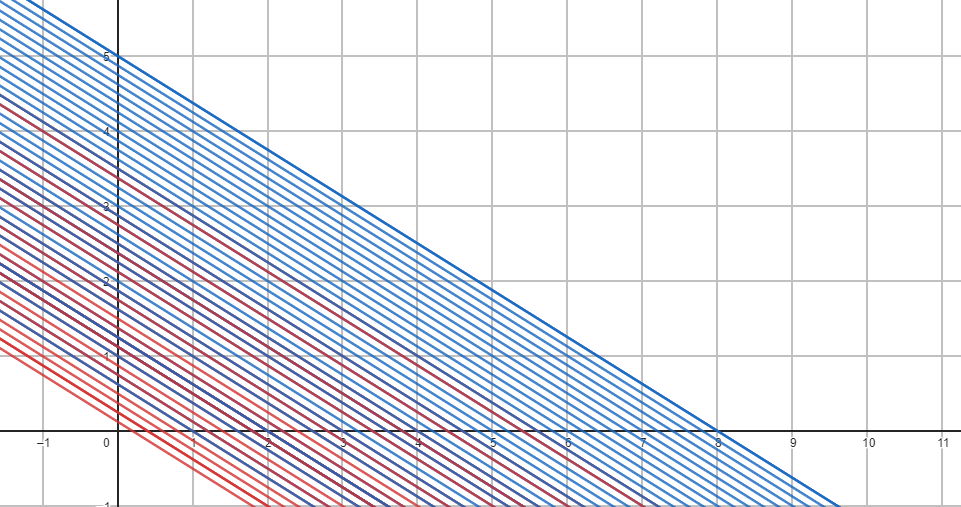}
\end{center}
\caption{{\sc GeoGebra}'s visualization of $5x+8y=d$ for $1\leq d\leq 40$.\label{f2}}
\end{figure}

Some additional observations that will help the students into understanding the problem more deeply are the following:
\begin{enumerate}
  \item For fixed $a$ and $b$, the various lines $\Li_d$ (where $d\in\Z_{\geq 0}$) form an infinite family of parallel lines.
  \item As $d$ grows, the line segment, where this line intersects with the first quadrant and the axis, also grows in length. This means that, intuitively, for higher values of $d$, the corresponding line should contain more desired points. Although this argument is not entirely strict (nor true) it still point towards the existence of $\f_{a,b}$.
  \item As $a$ and $b$ grow larger, the line segment, where $\Li_d$ intersects with the first quadrant and the axis, shrinks in length. This means, intuitively, that it should be more difficult to have common points with $\Lambda$. As a result the number $\f_{a,b}$ (if it exists) should probably be an increasing function of both $a$ and $b$.
\end{enumerate}
As a verification of these observations, it may be helpful to demonstrate these facts with explicit examples, or rather encourage students to experiment on their own and make their own conjectures, with the help of a computer program such as {\sc GeoGebra}. For example, see Figures~\ref{f1} and \ref{f2}.

Next, we proceed with the existence of $\f_{a,b}$. Following Proposition~\ref{gen_sol}, if $(x_0,y_0)$ is a solution of Eq.~\eqref{dioph}, then $d\in\C_{a,b}$ if, for some $k$, the point $(x_0-kb,y_0+ka)$ is in the first quadrant, or on some axis. Since this point cannot lie in the fourth quadrant, this is equivalent with
\[ (x_0-kb)(y_0+ka) \geq 0 \iff -ab\left( k+\frac{y_0}{a} \right)\left( k-\frac{x_0}{b} \right) \geq 0, \]
for some $k\in\Z$. The latter immediately yields an easy characterization of $\C_{a,b}$. In particular, $d\in\C_{a,b}$ if and only if there exist some integer $k$ in the interval $[-y_0/a , x_0/b]$ (or the interval $[x_0/b, -y_0/a]$ if $-y_0/a>x_0/b$), where $(x_0,y_0)$ is a solution of Eq.~\eqref{dioph}. Further, we observe that this condition is always satisfied when
\[ \left| \frac{y_0}{a} + \frac{x_0}{b} \right| \geq 1 \stackrel{ax_0+by_0=d}\Longleftrightarrow d\geq ab . \]
We have, so far, obtained a characterization of $\C_{a,b}$ and the facts that $\f_{a,b}$ exists and that $\f_{a,b}<ab$. Having obtained these, we are left with the computation of $\f_{a,b}$.

Towards this end, it may be useful to employ computers as a mean to visualize the problem, using various appropriate choices of $a$ and $b$, taking all $d\leq ab$, since the case $d>ad$ is settled, see Figures~\ref{f1} and \ref{f2}. In particular, the instructor is advised to proceed to the following key observations:
\begin{enumerate}
  \item The line $\Li_{ab}$ does not have common points with $\Lambda$ in the first quadrant, but one with each axis, namely $A(0,a)$ and $B(b,0)$. This is fairly easy to also explain theoretically.
  \item The lines right below $\Li_{ab}$ appear to have exactly one point in common with $\Lambda$ in the first quadrant.
  \item The first line that coincides with $\Lambda$ only in the second and the third quadrant seems to always be $\Li_{ab-a-b}$. In addition, among the points that this line shares with $\Lambda$ are $C(b-1,-1)$ and $D(-1,a-1)$. These facts are, again, fairly easy to theoretically explain.
\end{enumerate}
The above observations, not only suggest Eq.~\eqref{f}, but also suggest that the following classic geometric result could prove useful.
\begin{theorem}[Pick's theorem] \label{pick}
  Suppose that a polygon has its vertices in $\Lambda$. If $\mathcal I$ is the number of points of $\Lambda$ interior to the polygon and $\mathcal B$ is the number of points of $\Lambda$ on its boundary and $\mathcal A$ is the area of the polygon, then
  \begin{equation}\label{eq:pick} \mathcal A = \mathcal I + \frac{\mathcal B}2 - 1 . \end{equation}
\end{theorem}
The teacher is advised to bring Theorem~\ref{pick} to the students' attention at this point. 
It is natural to apply Theorem~\ref{pick} on the parallelogram $ABCD$, then, using the notation of the statement of Theorem~\ref{pick}, we easily deduce that $\mathcal B=4$ and $\mathcal A = a+b$, thus Eq.~\eqref{eq:pick} implies $\mathcal I = a+b-1$, in other words, the number of points of $\Lambda$ that are interior to $ABCD$ coincide with the number of lines $\Li_d$ (for $ab-a-b< d< ab$) that run through $ABCD$. Additionally, note that all these points lie within the first quadrant and that the point $(x,y)$ belongs to the line $\Li_{ax+by}$ (with $ab+by\in\Z_{\geq 0}$). Next, assume that the line $\Li_d$, with $ab-(a+b)<d<ab$, contains two interior lattice points of $ABCD$, say $(x_1,y_1)$ and $(x_2,y_2)$. Then
\[ \left. \begin{array}{r}
  d=ax_1+by_1 \\ d = ax_2+by_2
\end{array} \right\} \Rightarrow a(x_1-x_2) = b(y_2-y_1) \stackrel{\gcd(a,b)=1}{\Rightarrow} b\mid x_1-x_2 \text{ and } a\mid y_2-y_1 . \]
However, given that $0<x_1,x_2<b$ and $0<y_1,y_2<a$ the above is possible only if $(x_1,y_1)=(x_2,y_2)$, a contradiction. It follows that every line $\Li_d$, with $ab-(a+b)<d<ab$, contains at most one internal lattice point of $ABCD$. All the above facts, combined with the pigeonhole principle, imply that there is an one-to-one correspondence between internal lattice points of $ABCD$ and the lines $\Li_d$, with $ab-(a+b)<d<ab$. Eq.~\eqref{f} follows.
\subsection{Inductive approach}
Last, but not least, an inductive approach could be a welcome addition. For a formal presentation of the proof that we outline here, we refer the reader to the authors' recent work \cite{kapetanakisrizos23b}. In this approach, the students can be encouraged to write a computer script that could, given appropriate $a$ and $b$, check whether $d\in\C_{a,b}$ for $d$ within a reasonable range, or perform this experiment by hand. For example, see Tables~\ref{t1} and \ref{t2}, that we compiled using the {\sc GeoGebra} software. This should lead to a safe guess of the number $\f_{a,b}$ for this particular choice of $a$ and $b$.
\begin{table}[h]
\begin{center}
  \begin{tabular}{ccccc}
  \hline
$x$ & $y$ & $d$ & $3x+7y=d$                 & Notes \\ \hline
 -  &  -  & $1$ &     -                     & Impossible \\
 -  &  -  & $2$ &     -     & Impossible \\
$1$ & $0$ & $3$ & $3\cdot 1+ 7\cdot 0 = 3$ & - \\
 -  &  -  & $4$ &     -     & Impossible \\
 -  &  -  & $5$ &     -     & Impossible \\
$2$ & $0$ & $6$ & $3\cdot 2+ 7\cdot 0 = 6$ & - \\
$0$ & $1$ & $7$ & $3\cdot 0+ 7\cdot 1 = 7$ & - \\
 -  &  -  & $8$ &     -     & Impossible \\
$3$ & $0$ & $9$ & $3\cdot 3+ 7\cdot 0 = 9$ & - \\
$1$ & $1$ & $10$ & $3\cdot 1+ 7\cdot 1 = 10$ & - \\
 -  &  -  & $11$ &    -      & Impossible \\
$4$ & $0$ & $12$ & $3\cdot 4+ 7\cdot 0 = 12$ & - \\
$2$ & $1$ & $13$ & $3\cdot 2+ 7\cdot 1 = 13$ & - \\
$0$ & $2$ & $14$ & $3\cdot 0+ 7\cdot 2 = 14$ & - \\
$5$ & $0$ & $15$ & $3\cdot 5+ 7\cdot 0 = 15$ & - \\
$3$ & $1$ & $16$ & $3\cdot 3+ 7\cdot 1 = 16$ & - \\
$1$ & $2$ & $17$ & $3\cdot 1+ 7\cdot 2 = 17$ & - \\
$6$ & $0$ & $18$ & $3\cdot 6+ 7\cdot 0 = 18$ & - \\
$4$ & $1$ & $19$ & $3\cdot 4+ 7\cdot 1 = 19$ & - \\
$2$ & $2$ & $20$ & $3\cdot 2+ 7\cdot 2 = 20$ & - \\
$7$ & $0$ & $21$ & $3\cdot 7+ 7\cdot 0 = 21$ & - \\
$0$ & $3$ & $21$ & $3\cdot 0+ 7\cdot 3 = 21$ & - \\
$5$ &	$1$	& $22$ &	$3\cdot 5+7\cdot 1=22$ & - \\
$3$ & $2$	& $23$ & $3\cdot 3+7\cdot 2=23$ & - \\
$1$ & $3$ & $24$ & $3\cdot 1+7\cdot 3=24$ & - \\
$8$	& $0$	& $24$ & $3\cdot 8+7\cdot 0=24$ & - \\
$6$	& $1$	& $25$ & $3\cdot 6+7\cdot 1=25$ & - \\
$4$	& $2$	& $26$ &	$3\cdot 4+7\cdot 2=26$ & - \\
$2$ &	$3$	& $27$ &	$3\cdot 2+7\cdot 3=27$ & - \\
$9$	 & $0$ &	$27$	& $3\cdot 9+7\cdot 0=27$ & - \\
$7$ &	$1$ &	$28$	& $3\cdot 7+7\cdot 1=28$ & - \\
$0$	& $4$	& $28$	& $3\cdot 0+7\cdot4=28$ & - \\
$5$ & $2$ &	$29$	& $3\cdot 5+7\cdot 2=29$ & - \\
$3$ & 	$3$ & $30$ &	$3\cdot 3+7\cdot 3=30$ & - \\
$10$ &	$0$ &	$30$ &	$3\cdot 10+7\cdot 0=30$ & -
  \end{tabular}
\end{center}
\caption{Non-negative solutions of $3x+7y=d$, for $1\leq d\leq 30$.\label{t1}}
\end{table}

\begin{table}[h]
\begin{center}
  \begin{tabular}{ccccc}
  \hline
$x$ & $y$ & $d$ & $5x+8y=d$                 & Notes \\ \hline
 -  &  -  & $1$ &     -                     & Impossible \\
 -  &  -  & $2$ &     -     & Impossible \\
 -  &  -  & $3$ &     -     & Impossible \\
 -  &  -  & $4$ &     -     & Impossible \\
$1$ & $0$ & $5$ & $5\cdot 1+ 8\cdot 0 = 5$ & - \\
 -  &  -  & $6$ &     -     & Impossible \\
 -  &  -  & $7$ &     -     & Impossible \\
$0$ & $1$ & $8$ & $5\cdot 0+ 8\cdot 1 = 8$ & - \\
 -  &  -  & $9$ &     -     & Impossible \\
$2$ & $0$ & $10$ & $5\cdot 2+ 8\cdot 0 = 10$ & - \\
 -  &  -  & $11$ &     -     & Impossible \\
 -  &  -  & $12$ &     -     & Impossible \\
$1$ & $1$ & $13$ & $5\cdot 1+ 8\cdot 1 = 13$ & - \\
 -  &  -  & $14$ &     -     & Impossible \\
$3$ & $0$ & $15$ & $5\cdot 3+ 8\cdot 0 = 15$ & - \\
$0$ & $2$ & $16$ & $5\cdot 0+ 8\cdot 2 = 16$ & - \\
 -  &  -  & $17$ &     -     & Impossible \\
$2$ & $1$ & $18$ & $5\cdot 2+ 8\cdot 1 = 18$ & - \\
 -  &  -  & $19$ &     -     & Impossible \\
$4$ & $0$ & $20$ & $5\cdot 4+ 8\cdot 0 = 20$ & - \\
$1$ & $2$ & $21$ & $5\cdot 1+ 8\cdot 2 = 21$ & - \\
 -  &  -  & $22$ &     -     & Impossible \\
$3$ & $1$ & $23$ & $5\cdot 3+ 8\cdot 1 = 23$ & - \\
$0$ & $3$ & $24$ & $5\cdot 0+ 8\cdot 3 = 24$ & - \\
$5$ & $0$ & $25$ & $5\cdot 5+ 8\cdot 0 = 25$ & - \\
$2$ & $2$ & $26$ & $5\cdot 2+ 8\cdot 2 = 26$ & - \\
 -  &  -  & $27$ &     -     & Impossible \\
$4$ & $1$ & $28$ & $5\cdot 4+ 8\cdot 1 = 28$ & - \\
$1$ & $3$ & $29$ & $5\cdot 1+ 8\cdot 3 = 29$ & - \\
$6$ & $0$ & $30$ & $5\cdot 6+ 8\cdot 0 = 30$ & - \\
$3$ & $2$ & $31$ & $5\cdot 3+ 8\cdot 2 = 31$ & - \\
$0$ & $4$ & $32$ & $5\cdot 0+ 8\cdot 4 = 32$ & - \\
$5$ & $1$ & $33$ & $5\cdot 5+ 8\cdot 1 = 33$ & - \\
$2$ & $3$ & $34$ & $5\cdot 2+ 8\cdot 3 = 34$ & - \\
$7$ & $0$ & $35$ & $5\cdot 7+ 8\cdot 0 = 35$ & - \\
$4$ & $2$ & $36$ & $5\cdot 4+ 8\cdot 2 = 36$ & - \\
$1$ & $4$ & $37$ & $5\cdot 1+ 8\cdot 4 = 37$ & - \\
$6$ & $1$ & $38$ & $5\cdot 6+ 8\cdot 1 = 38$ & - \\
$3$ & $3$ & $39$ & $5\cdot 3+ 8\cdot 3 = 39$ & - \\
$8$ & $0$ & $40$ & $5\cdot 8+ 8\cdot 0 = 40$ & - \\
$0$ & $5$ & $40$ & $5\cdot 0+ 8\cdot 5 = 40$ & - 
  \end{tabular}
\end{center}
\caption{Non-negative solutions of $5x+8y=d$ for $1\leq d\leq 40$.\label{t2}}
\end{table}

After repeating the aforementioned procedure for a respectable number of choices for $a$ and $b$, the students may be able to notice a pattern, perhaps with some hints from the teacher/instructor with questions of the form:
\begin{itemize}
  \item[-] \emph{Does $\f_{a,b}$ appear to be larger or smaller than $ab$?}
  \item[-] \emph{Does $\f_{a,b}$ appear to be larger or smaller than $ab-a$?}
  \item[-] \emph{Does $\f_{a,b}$ appear to be larger or smaller than $ab-2a$?}
  \item[-] \ldots
\end{itemize}
Eventually, the students should recognize the pattern and guess Eq.~\eqref{f}. Then, proving Eq.~\eqref{f} can be an intriguing exercise. It is suggested to first observe that Proposition~\ref{gen_sol} implies the existence of two ``minimal'' pairs $(x_1,y_1)$ and $(x_2,y_2)$, with $x_1,y_2>0$ and $y_1,x_2<0$ and $ax_1+by_1=ax_2+by_2=1$. These expressions of $1$ are special, so we will call them \emph{minimal unit expressions}.

The idea is that, given an expression 
\begin{equation} \label{eq:d}
d=ax_0+by_0,
\end{equation}
with $x_0,y_0\geq 0$, a promising strategy, in order to obtain a nice expression for $d+1$ would be adding a minimal unit expression to Eq.~\eqref{eq:d}, as these expressions seem to be those that will cause the ``least'' damage to the non-negativity of $x_0,y_0$. In other words, at least one of these expressions should preserve their precious non-negativity, of course, given that $d$ is large enough.

With this overall strategy in mind, the instructor can proceed to the technical part. It is advised to split it into four parts at this point:
\begin{enumerate}
  \item Show that $ab-a-b\not\in\C_{a,b}$.
  \item Show that $ab-a-b+1\in\C_{a,b}$, by adding the suitable minimal unit expression to $ab-a-b = a(b-1) + b(-1)$. Notice that, since $y_1<0$, the only option is, in fact, to add the second minimal unit expression, that is,
  \[ ab-a-b+1 = a(b-1+x_2) + b(-1+y_2) . \]
  The students should, after studying Proposition~\ref{gen_sol}, be in position to understand why both factors of $a$ and $b$ above are non-negative.
  \item Show that one of the expressions
  \[ ab-a-b+2 = a(b-1+x_2+x_i) + b(-1+y_2+y_i), \]
  for $i=1,2$, has non-negative coefficients for both $a$ and $b$. This should provide some insight into how the induction step works.
  \item Show, using induction, that if 
  \[ k=a(x_0+x_i) + b(y_0+y_i) , \]
  where $i=1$ or $2$ and both coefficients of $a$ and $b$ are non-negative, then
  \[ k+1=a(x_0+x_i+x_j) + b(y_0+y_i+y_j) \]
  is such that  both coefficients of $a$ and $b$ are non-negative, for $j=1$ or $j=2$, for any $k\geq ab-a-b+2$. The induction should be on $k$ and the first step is already proven above.
\end{enumerate}
The above steps imply Eq.~\eqref{f}.
\section{Further work}
Should the students show enough interest on the Frobenius coin problem and one of the aforementioned approaches, the instructor can try to harden the problem further, by asking the students the following question:
\begin{quote}
  \emph{For fixed, relatively prime, $a,b\in\Z_{>0}$ and $d\in\C_{a,b}$, how many pairs $(x,y)\in\Z_{\geq 0}^2$ are there such that $d=ax+by$?}
\end{quote}
Equivalently, in how many ways can we obtain $d$ momentary units with coins of values $a$ and $b$ (with the usual restrictions on $a$ and $b$)? A closer observation to the aforementioned approaches, should provide the students enough insight to tackle this extension of the Frobenius coin problem and even solve it!
\section{Conclusive remarks}
Basic Number Theory, that is the theoretical preoccupation with the properties of integers and the equivalence classes of modulo, has been practically eliminated from secondary education (at least in Greek education system) and does not seem to attract the interest of researchers in Mathematics Education, while at the same time the relevant knowledge is used only as a prerequisite for mathematical competitions and Mathematics Olympiads. As much of the teaching research seems to have been aimed, in recent years, at new ways of formalizing and evaluating cognitive concepts, one explanation for this lack of interest could be that Number Theory is not easily formalized, as it moves simultaneously in the logic of discovery and that of proof, and its reasoning is largely based on common sense, which is not easily standardize. However, taking into consideration our recent experience \cite{rizosgkrekas23}, we believe that students are inclined to think intuitively, discover patterns, and give reasons on problems in basic Number Theory, such as those suggested above. A future study based on a teaching intervention utilizing tasks such as those presented in Section~\ref{motivation}, could explore this possibility in detail.
	\bibliographystyle{abbrv}
	\bibliography{references}
\end{document}